\documentclass[12pt, reqno]{amsart}
\usepackage{amsmath, amstext, amsbsy, amssymb, amscd}

\newtheorem{theorem}{Theorem}[section]
\newtheorem{lemma}[theorem]{Lemma}
\newtheorem{proposition}[theorem]{Proposition}
\newtheorem{corollary}[theorem]{Corollary}
\theoremstyle{definition}
\newtheorem{definition}[theorem]{Definition}

\theoremstyle{remark}
\newtheorem{remark}[theorem]{Remark}
\newtheorem{assumption}[theorem]{Assumption}
\newtheorem{conjecture}[theorem]{Conjecture}

\numberwithin{equation}{section}

\newcommand{\C}{\mathbb C}

\newcommand{\Hom}{{\rm Hom}}

\newcommand{\Pee}{\mathbb P}

\newcommand{\Supp}{{\rm Supp}}

\newcommand{\w}{\tilde}

\newcommand{\Xn}{X^{[n]}}
\newcommand{\Z}{\mathbb Z}

{\vskip-\lastskip\medskip
  \noindent
  {\em #1.}\enspace
  }%
{\qed\par\medskip
  }

\begin{document}
\title[Moduli spaces of ideal sheaves and Donaldson-Thomas invariants]
{On certain moduli spaces of ideal sheaves and
Donaldson-Thomas invariants}

\author[Sheldon Katz]{Sheldon Katz$^1$}
\address{Department of Mathematics, University of Illinois at
Urbana-Champaign, Urbana, IL 61801} \email{katz@math.uiuc.edu}
\thanks{${}^1$Partially supported by an NSF grant}

\author[Wei-Ping Li]{Wei-Ping Li$^2$}
\address{Department of Mathematics, HKUST, Clear Water Bay,
Kowloon, Hong Kong} \email{mawpli@ust.hk}
\thanks{${}^2$Partially supported by the grant HKUST6114/02P}

\author[Zhenbo Qin]{Zhenbo Qin$^3$}
\address{Department of Mathematics, University of Missouri,
Columbia, MO 65211, USA} \email{zq@math.missouri.edu}
\thanks{${}^3$Partially supported by an NSF grant}


\begin{abstract}
We determine the structure of certain moduli spaces of ideal
sheaves by generalizing an earlier result of the first author.
As applications, we compute the (virtual) Hodge polynomials
of these moduli space, and calculate the Donaldson-Thomas
invariants of certain $3$-folds with trivial canonical classes.
\end{abstract}

\maketitle
\date{}

\section{Introduction}
\label{sect_intro}

Donaldson-Thomas theory was introduced in \cite{DT, Tho} via
integrals over the moduli spaces of semistable sheaves and
via the theory of virtual fundamental cycles. It was further
developed by Maulik-Pandharipande \cite{MP} and Jun Li.
In \cite{LQ1}, the second and third authors
constructed rank-$2$ stable vector bundles
over certain Calabi-Yau manifolds, and calculated the
corresponding Donaldson-Thomas invariants.
Recently, a rich interplay among Donaldson-Thomas theory,
Gromov-Witten theory and Gopakumar-Vafa invariants
has been proposed, and in certain cases, verified in
\cite{MNOP1, MNOP2, Kat, OP} (see the references
there for other papers). For this interplay,
Donaldson-Thomas theory is defined via the moduli spaces
parametrizing ideal sheaves of $1$-dimensional closed
subschemes in a $3$-fold, while Gromov-Witten theory
is based on the moduli spaces of stable morphisms to
the same $3$-fold.  A complete mathematical theory of
Gopakumar-Vafa invariants has not yet been developed.

In this paper, we study the moduli spaces of ideal sheaves
in certain $3$-folds $X$, and compute the corresponding
Donaldson-Thomas invariants of $X$.
Let $S$ be a smooth projective surface, and $\mu: X \to S$
be a Zariski-locally trivial fibration whose fibers are
smooth irreducible curves of genus-$g$.
Let $\beta \in H_2(X; \Z)$ be the class of a fiber.
For two nonnegative integers $m$ and $n$, the moduli space
$$\mathfrak I_{m(1-g)+n}(X, m\beta)$$
parametrizes ideal sheaves $I_Z \subset \mathcal O_X$
where $Z$ denotes $1$-dimensional closed
subschemes of $X$ with $\chi(\mathcal O_Z)=m(1-g)+n$ and
$[Z] = m\beta$. Here $[Z]$ is the class associated to
the dimension-$1$ component (weighted by their
intrinsic multiplicities) of $Z$. When $n = 0$,
the moduli space $\mathfrak I_{m(1-g)}(X, m\beta)$ is
naturally identified with the Hilbert scheme $S^{[m]}$
which parametrizes the length-$m$ $0$-dimensional
closed subschemes of $S$. Our main result is to
determine $\mathfrak I_{m(1-g)+1}(X, m\beta)$.

\begin{theorem}   \label{intro_thm}
Identify $S^{[m]}$ and $\mathfrak I_{m(1-g)}(X, m\beta)$.
Let $\widetilde{S^{[m]} \times X}$ be the blow-up of
$S^{[m]} \times X$ along the universal curve. Then,
we have an isomorphism
\begin{eqnarray*}
\mathfrak I_{m(1-g)+1}(X, m\beta)
\cong \widetilde{S^{[m]} \times X}.
\end{eqnarray*}
Moreover, the moduli space
$\mathfrak I_{m(1-g)+1}(X, m\beta)$ is smooth.
\end{theorem}

As an application to Theorem~\ref{intro_thm}, we compute
the (virtual) Hodge polynomials of the moduli space
$\mathfrak I_{m(1-g)+1}(X, m\beta)$ using the results
of Cheah and Ellingsrud-Str\o mme \cite{Che, ES}.
As another application, we show that when $X = S \times C$
with $K_S = 0$ and $K_C = 0$, the Donaldson-Thomas invariant
of the $3$-fold $X$ corresponding to the moduli space
$\mathfrak I_{m(1-g)+1}(X, m\beta)$ is equal to zero.

Theorem~\ref{intro_thm} is proved in Sect.~\ref{sect_moduli}.
The key step is to establish Lemma~\ref{01} which provides
a relation between the two moduli spaces $\mathfrak I_{1-g}(X, \beta)$
and $\mathfrak I_{2-g}(X, \beta)$. By assuming a canonical
isomorphism between the blow-up along the universal curve
and the projectivization associated to the ideal sheaf
defining the universal curve, our Lemma~\ref{01}
generalizes an earlier result of the first author
(see Remark~\ref{katz_rmk}). Indeed,
the proof of Lemma~\ref{01} follows the same argument as
in the proof of the Lemma~1 in \cite{Kat}.
The only difference is that instead of using local arguments,
we apply the universal properties of various constructions.

In Sect.~\ref{sect_com}, we discuss some relation between
the two moduli spaces $\mathfrak I_{1-g}(X, \beta)$
and $\mathfrak I_{3-g}(X, \beta)$. We show that
$\mathfrak I_{3-g}(X, \beta)$ is not smooth in general.

\bigskip\noindent
{\bf Acknowledgments.} The authors thank Jun Li for
valuable discussions.
\section{\bf The moduli space $\mathfrak I_{2-g}(X, \beta)$
   and main results}
\label{sect_moduli}

\subsection{\bf The definition of the moduli space
                $\mathfrak I_n(X, \beta)$}
\label{subsect_def}
$\,$
\smallskip

Let $X$ be a smooth projective complex variety of dimension $r$.
For a $1$-dimensional closed subscheme $Z$ of $X$, let $[Z]
\in H_2(X; \Z)$ be the class associated to the dimension-$1$
component (weighted by their intrinsic multiplicities) of $Z$.
As in \cite{MNOP1, MNOP2}, we make the following definition.

\begin{definition} \label{def_moduli}
For a fixed homology class $\beta \in H_2(X; \Z)$ and a fixed integer
$n$, we define $\mathfrak I_n(X, \beta)$ to be the moduli space
parametrizing the ideal sheaves $I_Z$ of $1$-dimensional\footnote{If
$\beta=0$, the closed subschemes $Z$ are actually $0$-dimensional.
To avoid repeatedly having a separate discussion of this case, we abuse
terminology slightly in this paper by speaking only of $1$-dimensional
closed schemes.} closed subschemes $Z$ of
$X$ satisfying the conditions:
\begin{eqnarray}   \label{In}
\chi(\mathcal O_Z)=n, \qquad [Z] = \beta.
\end{eqnarray}
\end{definition}

Notice that it is also convenient to regard
$\mathfrak I_n(X, \beta)$ as the moduli space parametrizing
the corresponding closed subschemes $Z$.

The degree-$0$ moduli space $\mathfrak I_n(X, 0)$ is isomorphic
to the Hilbert scheme $\Xn$ parametrizing length-$n$
$0$-dimensional closed subschemes of $X$. In general,
when $\beta \ne 0$, the moduli space $\mathfrak I_n(X, \beta)$
is only part of the Hilbert scheme defined in terms of certain
degree-$1$ Hilbert polynomial (see \cite{Gro}).
The Zariski tangent space of $\mathfrak I_n(X, \beta)$
at a point $[Z]$ is canonically isomorphic to
\begin{eqnarray}   \label{Zar_tang}
\text{Hom}(I_Z, \mathcal O_Z).
\end{eqnarray}
By the Lemma 1 in \cite{MNOP2}, when $\dim (X) = 3$,
the perfect obstruction theory on $\mathfrak I_n(X, \beta)$
defined in \cite{Tho} has virtual dimension equal to
\begin{eqnarray}   \label{dim_In}
-(\beta \cdot K_X)
\end{eqnarray}
where $K_X$ stands for the canonical class of $X$.

\subsection{\bf A relation between $\mathfrak I_{1-g}(X, \beta)$
and $\mathfrak I_{2-g}(X, \beta)$}
\label{subsect_rel}
$\,$
\smallskip

Let $g \ge 0$. Assume that all the closed subschemes parametrized
by $\mathfrak I_{1-g}(X, \beta)$ are curves of arithmetic genus $g$,
and that the $1$-dimensional components of all the closed
subschemes parametrized by $\mathfrak I_{2-g}(X, \beta)$
are supported on the curves parametrized by
$\mathfrak I_{1-g}(X, \beta)$. It follows that if
$[Z] \in \mathfrak I_{2-g}(X, \beta)$, then there exists
a unique $[Z'] \in \mathfrak I_{1-g}(X, \beta)$
together with an exact sequence:
\begin{eqnarray}   \label{exact}
0 \to I_Z \to I_{Z'} \to \mathcal O_x \to 0
\end{eqnarray}
where $x$ is some point in $X$. Let $\mathcal I_{1-g}$ be
the universal ideal sheaf over
$$\mathfrak I_{1-g}(X, \beta) \times X,$$
and let $\mathcal C_{1-g} \subset \mathfrak I_{1-g} \times X$
be the universal curve. Then, $\mathcal I_{1-g} =
I_{\mathcal C_{1-g}}$.

Let $\Pee(\mathcal I_{1-g})$ be the projectivization of
the sheaf $\mathcal I_{1-g}$, and
\begin{eqnarray}   \label{til_pi}
\tilde \pi: \Pee(\mathcal I_{1-g})
\to \mathfrak I_{1-g}(X, \beta) \times X
\end{eqnarray}
be the natural projection. Then there exists a universal quotient
\begin{eqnarray}   \label{P_univ}
\tilde{\pi}^*\mathcal I_{1-g} \to
\mathcal O_{\Pee(\mathcal I_{1-g})}(1) \to 0
\end{eqnarray}
over $\Pee(\mathcal I_{1-g})$. Let
\begin{eqnarray}   \label{pi}
\pi: \mathfrak I_{1-g}(X \widetilde{, \beta) \times X}
\to \mathfrak I_{1-g}(X, \beta) \times X
\end{eqnarray}
be the blow-up of $\mathfrak I_{1-g}(X, \beta) \times X$
along $\mathcal C_{1-g}$, and $E$ be the exceptional divisor.
Then there exists a surjection over
$\mathfrak I_{1-g}(X \widetilde{, \beta) \times X}$:
\begin{eqnarray}   \label{pi*}
\pi^*\mathcal I_{1-g} \to \mathcal O_{\mathfrak I_{1-g}
  (X \widetilde{, \beta) \times X}}(-E) \to 0.
\end{eqnarray}
By the universal property of the projectivizations,
there exists a canonical morphism
$\phi: \mathfrak I_{1-g}(X \widetilde{, \beta) \times X}
\to \Pee(\mathcal I_{1-g})$ making a commutative diagram:
\begin{eqnarray}   \label{blowup_grass.1}
\begin{array}{ccccc}
\mathfrak I_{1-g}(X \widetilde{, \beta) \times X}
  &&\overset \phi \longrightarrow&&\Pee(\mathcal I_{1-g})\\
&\searrow \pi&&\tilde{\pi}\swarrow&\\
&&\mathfrak I_{1-g}(X, \beta) \times X&&\\
\end{array}
\end{eqnarray}
such that the pull-back of (\ref{P_univ}) via $\phi^*$
is (\ref{pi*}). In particular,
\begin{eqnarray}   \label{blowup_grass.2}
\phi^*\mathcal O_{\Pee(\mathcal I_{1-g})}(1)
= \mathcal O_{\mathfrak I_{1-g}
  (X \widetilde{, \beta) \times X}}(-E).
\end{eqnarray}

\begin{lemma}  \label{01}
Assume that the closed subschemes parametrized
by $\mathfrak I_{1-g}(X, \beta)$ are curves of arithmetic genus $g$,
and that the $1$-dimensional components of the closed
subschemes parametrized by $\mathfrak I_{2-g}(X, \beta)$
are supported on the curves parametrized by
$\mathfrak I_{1-g}(X, \beta)$. If the canonical morphism $\phi$
in (\ref{blowup_grass.1}) is an isomorphism, then
\begin{eqnarray}   \label{01.0}
\mathfrak I_{2-g}(X, \beta) \cong
\mathfrak I_{1-g}(X \widetilde{, \beta) \times X}.
\end{eqnarray}
\end{lemma}
\begin{proof}
For convenience, let $\mathfrak I_n = \mathfrak I_n(X, \beta)$.
For $i = 2$ or $3$, let
$$\pi_{1i}: \mathfrak I_{1-g} \times X \times X
\to \mathfrak I_{1-g} \times X$$
be the projection to the first and $i$-th factors.
Let $\Delta_X$ be the diagonal
in $X \times X$. Over $\mathfrak I_{1-g} \times X \times X$,
we have the composition of the natural morphisms:
\begin{eqnarray}   \label{01.1}
\pi_{13}^*\mathcal I_{1-g} \to
\mathcal O_{\mathfrak I_{1-g} \times X \times X} \to
\mathcal O_{\mathfrak I_{1-g} \times \Delta_X}.
\end{eqnarray}
Regarding (\ref{01.1}) as a family of
morphisms parametrized by $\mathfrak I_{1-g} \times X$,
we see that (\ref{01.1}) vanishes precisely along
$\mathcal C_{1-g}$. So there is an induced morphism:
\begin{eqnarray}   \label{01.2}
(\pi \times \text{Id}_X)^*\pi_{13}^*\mathcal I_{1-g} \to
(\pi \times \text{Id}_X)^*\mathcal O_{\mathfrak I_{1-g}
\times \Delta_X}(-\w \pi_{12}^*E)
\end{eqnarray}
over $\widetilde{\mathfrak I_{1-g} \times X} \times X$,
where $\w \pi_{12}: \widetilde{\mathfrak I_{1-g} \times X}
\times X \to \widetilde{\mathfrak I_{1-g} \times X}$ is
the natural projection. Note that we have a commutative
diagram of morphisms:
\begin{eqnarray}   \label{add1}
\begin{array}{ccc}
\widetilde{\mathfrak I_{1-g} \times X} \times X
  &\overset {\w \pi_{12}} \longrightarrow
  &\widetilde{\mathfrak I_{1-g} \times X}\\
\qquad\qquad\quad \downarrow \pi \times \text{Id}_X
  &&\downarrow \pi\\
\mathfrak I_{1-g} \times X \times X
  &\overset {\pi_{12}} \longrightarrow
  &\mathfrak I_{1-g} \times X.
\end{array}
\end{eqnarray}
Therefore we see from (\ref{pi*})
that (\ref{01.2}) is surjective since
\begin{eqnarray*}
& &(\pi \times \text{Id}_X)^*\pi_{13}^*\mathcal I_{1-g}
   |_{(\pi \times \text{Id}_X)^{-1}(\mathfrak I_{1-g}
   \times \Delta_X)}  \\
&=&(\pi \times \text{Id}_X)^*\pi_{12}^*\mathcal I_{1-g}
   |_{(\pi \times \text{Id}_X)^{-1}(\mathfrak I_{1-g}
   \times \Delta_X)}  \\
&=&\w \pi_{12}^* \pi^*\mathcal I_{1-g}
   |_{(\pi \times \text{Id}_X)^{-1}(\mathfrak I_{1-g}
   \times \Delta_X)}  \\
&\twoheadrightarrow&\w \pi_{12}^*
   \mathcal O_{\widetilde{\mathfrak I_{1-g} \times X}}(-E)
   |_{(\pi \times \text{Id}_X)^{-1}
   (\mathfrak I_{1-g} \times \Delta_X)} \\
&=&(\pi \times \text{Id}_X)^*\mathcal O_{\mathfrak I_{1-g}
   \times \Delta_X}(-\w \pi_{12}^*E).
\end{eqnarray*}
Let $\tilde{\mathcal I}$ be the kernel of (\ref{01.2}).
Then we have an exact sequence
\begin{eqnarray}   \label{01.3}
0 \to \tilde{\mathcal I} \to
(\pi \times \text{Id}_X)^*\pi_{13}^*\mathcal I_{1-g} \to
(\pi \times \text{Id}_X)^*\mathcal O_{\mathfrak I_{1-g}
\times \Delta_X}(-\w \pi_{12}^*E) \to 0
\end{eqnarray}
over $\widetilde{\mathfrak I_{1-g} \times X} \times X$.
Now $\tilde{\mathcal I}$ is flat over
$\widetilde{\mathfrak I_{1-g} \times X}$ since the other
two terms in (\ref{01.3}) are. Also,
the fibers of $\tilde{\mathcal I}$ over
$\widetilde{\mathfrak I_{1-g} \times X}$ are ideal sheaves
parametrized by $\mathfrak I_{2-g}$.

Next, we check that $\tilde{\mathcal I} \subset
\mathcal O_{\widetilde{\mathfrak I_{1-g} \times X} \times X}$
is the universal ideal sheaf. Let $\mathcal J \subset
\mathcal O_{T \times X}$ be a flat family of ideal sheaves
in $\mathfrak I_{2-g}$ parametrized by $T$.
Let $\mathcal J'$ be the saturation of $\mathcal J \subset
\mathcal O_{T \times X}$ (see Definition~1.1.5 in \cite{HL}).
Then, $\mathcal J'$ is a flat family of ideal sheaves in
$\mathfrak I_{1-g}$, and fits in an exact sequence
\begin{eqnarray}   \label{01.4}
0 \to \mathcal J \to \mathcal J' \to
(\text{Id}_T \times f_1)_*\mathcal L \to 0
\end{eqnarray}
for some morphism $f_1: T \to X$ and some line bundle
$\mathcal L$ on $T$.
The flat family $\mathcal J'$ over $T \times X$ induces
a morphism $f_2: T \to \mathfrak I_{1-g}$ such that
\begin{eqnarray}   \label{01.5}
(f_2 \times \text{Id}_X)^*\mathcal I_{1-g} = \mathcal J'.
\end{eqnarray}
By base-change, we have the commutative diagram:
\begin{eqnarray}    \label{01.6}
\begin{array}{ccccc}
&&\Pee(\mathcal J')&\longrightarrow&
  \Pee(\mathcal I_{1-g}) \\
&&\downarrow&&\downarrow\\
T&\overset {\text{Id}_T \times f_1} \longrightarrow
  &T \times X&\overset {f_2 \times \text{Id}_X}
  \longrightarrow &\mathfrak I_{1-g} \times X.\\
\end{array}
\end{eqnarray}
By (\ref{01.4}), we get a surjection
$(\text{Id}_T \times f_1)^*\mathcal J' \to
(\text{Id}_T \times f_1)^*(\text{Id}_T \times f_1)_*
\mathcal L \to 0$. Since the natural morphism
$(\text{Id}_T \times f_1)^*(\text{Id}_T \times f_1)_*
\mathcal L \to \mathcal L$ is surjective, we obtain
$$(\text{Id}_T \times f_1)^*\mathcal J' \to
\mathcal L \to 0$$
over $T$. By the universal property of $\Pee(\mathcal J')$,
we obtain a commutative diagram
\begin{eqnarray}    \label{01.7}
\begin{array}{ccccc}
&&\Pee(\mathcal J')&\longrightarrow&
  \Pee(\mathcal I_{1-g}) \\
&\nearrow&\downarrow&&\downarrow\\
T&\overset {\text{Id}_T \times f_1} \longrightarrow
  &T \times X&\overset {f_2 \times \text{Id}_X}
  \longrightarrow&\mathfrak I_{1-g} \times X.\\
\end{array}
\end{eqnarray}
Thus the morphism $f_2 \times f_1: T \to
\mathfrak I_{1-g} \times X$ can be lifted to a morphism
$$f: T \longrightarrow \Pee(\mathcal I_{1-g})
\overset {\phi^{-1}} \longrightarrow
\widetilde{\mathfrak I_{1-g} \times X}.$$

Finally, we apply the pull-back $(f \times \text{Id}_X)^*$
to (\ref{01.3}). Using (\ref{blowup_grass.2}) and
the property of the morphism $f$,
we obtain a commutative diagram:
\begin{eqnarray*}
\begin{array}{ccccccccc}
0&\to&(f \times \text{Id}_X)^*\tilde{\mathcal I}&\to&
  \mathcal J' &\to&(f \times \text{Id}_X)^*(\pi \times
  \text{Id}_X)^*\mathcal O_{\mathfrak I_{1-g}
  \times \Delta_X}(-\w \pi_{12}^*E)&\to&0  \\
&&&&\Vert&&\downarrow&        \\
0&\to&\mathcal J&\to&\mathcal J'&\to&
  (\text{Id}_T \times f_1)_*\mathcal L&\to&0
\end{array}
\end{eqnarray*}
over $T \times X$. In particular, there exists an injection
\begin{eqnarray}   \label{final}
0 \to (f \times \text{Id}_X)^*\tilde{\mathcal I}
\overset \psi \to \mathcal J.
\end{eqnarray}
Since both $(f \times \text{Id}_X)^*\tilde{\mathcal I}$
and $\mathcal J$ are flat families of ideal sheaves in
$\mathfrak I_{2-g}(X, \beta)$, we conclude that
the morphism $\psi$ is an isomorphism.
\end{proof}

\begin{remark}  \label{katz_rmk}
If the universal curve $\mathcal C_{1-g} \subset
\mathfrak I_{1-g} \times X$ is a local complete intersection,
then $(\mathcal I_{1-g})^n \cong
\text{\rm Sym}^n (\mathcal I_{1-g})$, and so the canonical
morphism $\phi$ in (\ref{blowup_grass.1}) is an isomorphism.
Hence Lemma~\ref{01} indeed generalizes the Lemma~1 in \cite{Kat}.
\end{remark}

\subsection{\bf Applications}
\label{subsect_app}
$\,$
\smallskip

In this subsection, we adopt the following basic
assumptions.

\begin{assumption} \label{assumption}
We assume that $X$ admits a Zariski-locally trivial fibration
\begin{eqnarray}   \label{fibration}
\mu: X \to S
\end{eqnarray}
where $S$ is a smooth surface, the fibers are smooth irreducible
curves of genus-$g$, and $\beta \in H_2(X; \Z)$ is
the class of a fiber. For $m, n \ge 0$, we let
\begin{eqnarray}   \label{Imn}
\mathfrak I_{m, n} := \mathfrak I_{m(1-g)+n}(X, m\beta).
\end{eqnarray}
\end{assumption}

Note that the elements in $\mathfrak I_{m, 0}$ correspond to the
ideal sheaves of the form $\mu^*I_\xi$ for some $\xi \in S^{[m]}$.
Hence there exists a bijective morphism $S^{[m]} \to
\mathfrak I_{m, 0}$. It is well-known that the Hilbert scheme
$S^{[m]}$ is smooth. Combining with (\ref{Zar_tang}),
one checks that the moduli space $\mathfrak I_{m, 0}$ is smooth
and that
\begin{eqnarray}   \label{Im0}
\mathfrak I_{m, 0} \cong S^{[m]}.
\end{eqnarray}

\begin{theorem}   \label{thm}
Identify $S^{[m]}$ and $\mathfrak I_{m, 0}$ by (\ref{Im0}).
Let $\widetilde{S^{[m]} \times X}$ be the blow-up of
$S^{[m]} \times X$ along the universal curve. Then,
we have an isomorphism
\begin{eqnarray*}
\mathfrak I_{m, 1} \cong \widetilde{S^{[m]} \times X}.
\end{eqnarray*}
Moreover, the moduli space $\mathfrak I_{m, 1}
= \mathfrak I_{m(1-g)+1}(X, m\beta)$ is smooth.
\end{theorem}
\begin{proof}
We will apply Lemma~\ref{01} to the present situation.
First of all, note that the $1$-dimensional components
of the closed subschemes parametrized by
$\mathfrak I_{m, 1}$ are supported on the curves
parametrized by $\mathfrak I_{m, 0}$. Next, let
$\mathcal Z_m \subset S^{[m]} \times S$ be the universal
codimension-$2$ subscheme. Set-theoretically,
\begin{eqnarray}   \label{thm.1}
\mathcal Z_m = \{(\xi, x) \in S^{[m]} \times S|\,\,
x \in \Supp(\xi) \}.
\end{eqnarray}
Let $\widetilde{S^{[m]} \times S}$ be the blow-up of
$S^{[m]} \times S$ along $\mathcal Z_m$.
By the results in \cite{ES}, $\widetilde{S^{[m]} \times S}$
is smooth and there exists a canonical isomorphism:
\begin{eqnarray}   \label{thm.2}
\widetilde{S^{[m]} \times S} \cong \Pee(I_{\mathcal Z_m}).
\end{eqnarray}

Now the universal curve $\mathcal C_m \subset S^{[m]}
\times X$ is $(\text{Id}_{S^{[m]}} \times \mu)^*
\mathcal Z_m$. Since $\mu$ is Zariski-locally trivial,
we obtain canonical isomorphisms:
\begin{eqnarray}
   \widetilde{S^{[m]} \times X}
&\cong&(\widetilde{S^{[m]} \times S})
   \times_{S^{[m]} \times S} (S^{[m]} \times X)
   \label{thm.3}  \\
&\cong&\Pee(I_{\mathcal Z_m})
   \times_{S^{[m]} \times S} (S^{[m]} \times X)
   \nonumber   \\
&\cong&\Pee(I_{\mathcal C_m}).   \nonumber
\end{eqnarray}
By Lemma~\ref{01}, we have an isomorphism
$\mathfrak I_{m, 1} \cong \widetilde{S^{[m]} \times X}$.
Note from the isomorphism (\ref{thm.3}) that
$\widetilde{S^{[m]} \times X}$ is smooth.
Hence $\mathfrak I_{m, 1}$ is smooth as well.
\end{proof}

\begin{corollary}   \label{cor_euler}
Let $e(\, \cdot \, ; s, t)$
denote the (virtual) Hodge polynomial. Then,
\begin{eqnarray*}
\sum_{m=0}^{+\infty} e(\mathfrak I_{m, 1}; s, t) q^m
= {q \over 1-stq} \cdot e(X; s, t) \cdot
\prod_{k=1}^{+\infty} \prod_{i, j} \left (
{1 \over 1-s^{i+k-1}t^{j+k-1}q^k} \right )^{e^{i, j}(S)}
\end{eqnarray*}
where $e^{i, j}(S) = (-1)^{i+j} h^{i, j}(S)$ and
$h^{i, j}(S)$ denotes the Hodge numbers of $S$.
\end{corollary}
\begin{proof}
Since $\mu$ is Zariski-locally trivial,
(\ref{thm.3}) implies that the natural projection
$$\widetilde{S^{[m]} \times X} \to
\widetilde{S^{[m]} \times S}$$
is a Zariski-locally trivial fibration with fibers
being isomorphic to the fibers $C$ of $\mu$. Thus,
$e(\widetilde{S^{[m]} \times X}; s, t) =
e(C; s, t) \cdot e(\widetilde{S^{[m]} \times S}; s, t)$.
By Theorem~\ref{thm},
\begin{eqnarray}   \label{cor_euler.1}
\sum_{m=0}^{+\infty} e(\mathfrak I_{m, 1}; s, t) q^m
= e(C; s, t) \cdot \sum_{m=0}^{+\infty}
  e(\widetilde{S^{[m]} \times S}; s, t) q^m.
\end{eqnarray}
By the Proposition~2.2 in \cite{ES}, $\widetilde{S^{[m]}
\times S}$ is isomorphic to the incidence variety:
\begin{eqnarray}   \label{cor_euler.2}
S_{m, m+1} = \{(\eta, \xi) \in S^{[m]} \times S^{[m+1]}|\,
\eta \subset \xi \}.
\end{eqnarray}
The (virtual) Hodge polynomial of $S_{m, m+1}$ has been
computed by Cheah:
\begin{eqnarray*}
\sum_{m=0}^{+\infty} e(S_{m, m+1}; s, t) q^m
= {q \over 1-stq} \cdot e(S; s, t) \cdot
\prod_{k=1}^{+\infty} \prod_{i, j} \left (
{1 \over 1-s^{i+k-1}t^{j+k-1}q^k} \right )^{e^{i, j}(S)}
\end{eqnarray*}
(see p.485 in \cite{Che}). Combining this with
(\ref{cor_euler.1}) completes the proof.
\end{proof}

\begin{remark}   \label{rmk}
It is well-known that $e(\, \cdot \, ; 1, 1)$ is equal
to the topological Euler number $\chi(\cdot)$.
Therefore, we see from Corollary~\ref{cor_euler} that
\begin{eqnarray}     \label{rmk.1}
\sum_{m=0}^{+\infty} \chi(\mathfrak I_{m, 1}) q^m
= {q \over 1 - q} \cdot \chi(X) \cdot
\prod_{k=1}^{+\infty} \left (
{1 \over 1-q^k} \right )^{\chi(S)}.
\end{eqnarray}
It is also interesting to note
from (\ref{Im0}) and G\" ottsche's formula in
\cite{Got} that
\begin{eqnarray}     \label{rmk.2}
\sum_{m=0}^{+\infty} \chi(\mathfrak I_{m, 0}) q^m
= \sum_{m=0}^{+\infty} \chi(S^{[m]}) q^m
= \prod_{k=1}^{+\infty} \left ( {1 \over 1-q^k}
  \right )^{\chi(S)}.
\end{eqnarray}
\end{remark}

Let $X = S \times C$ where $C$ is an elliptic curve and
$S$ is a smooth surface with $K_S = 0$. Let $\mu: X
\to S$ be the first projection. Then, $K_X = 0$.
By (\ref{dim_In}), the virtual dimension of the moduli
space $\mathfrak I_{m, n}$ is zero.
The corresponding {\it Donaldson-Thomas invariant} is an integer.
We denote this Donaldson-Thomas invariant by
\begin{eqnarray}     \label{dt}
N_{m, n}.
\end{eqnarray}

\begin{corollary}   \label{cor_dt}
Let $X = S \times C$ where $C$ is an elliptic curve and
$S$ is a smooth surface with $K_S = 0$. Then,
$N_{m, 1} = 0$ for every $m \ge 0$.
\end{corollary}
\begin{proof}
Since $K_X = 0$ and the moduli space $\mathfrak I_{m, 1}$
is smooth, the obstruction bundle over
$\mathfrak I_{m, 1}$ is the dual of the tangent bundle
of $\mathfrak I_{m, 1}$ (see \cite{Tho}). Hence
$$N_{m, 1} = (-1)^{\dim \mathfrak I_{m, 1}} \cdot
\chi(\mathfrak I_{m, 1}).$$
Since $\chi(X) = 0$, we see from (\ref{rmk.1}) that
$\chi(\mathfrak I_{m, 1}) = 0$. Therefore, $N_{m, 1} = 0$.
\end{proof}

\begin{conjecture}   \label{conj_dt}
Let $X = S \times C$ where $C$ is an elliptic curve and
$S$ is a smooth surface with $K_S = 0$. Then,
$N_{m, n} = 0$ for all $m \ge 0$ and $n \ge 1$.
\end{conjecture}

\section{\bf The moduli space $\mathfrak I_{3-g}(X, \beta)$}
\label{sect_com}

In this section, we make a few comments about the moduli
space $\mathfrak I_{3-g}(X, \beta)$. As in Subsect.~\ref{subsect_rel},
we assume that all the closed subschemes parametrized by
$\mathfrak I_{1-g}(X, \beta)$ are curves of arithmetic genus $g$,
and that the $1$-dimensional components of all the closed
subschemes parametrized by $\mathfrak I_{3-g}(X, \beta)$
are supported on the curves parametrized by
$\mathfrak I_{1-g}(X, \beta)$. Therefore, if $[Z] \in
\mathfrak I_{3-g}(X, \beta)$, then there exists a unique $[Z']
\in \mathfrak I_{1-g}(X, \beta)$ together with an exact sequence:
\begin{eqnarray}   \label{exact2}
0 \to I_Z \to I_{Z'} \to Q \to 0
\end{eqnarray}
where $Q$ is a torsion sheaf on $X$ with $\ell(Q) = 2$.

\begin{lemma} \label{Q_str}
Let $Q$ be a torsion sheaf on $X$ from (\ref{exact2}).
Then, either $Q \cong \mathcal O_{\xi}$ for some
$\xi\in X^{[2]}$, or $Q \cong \mathcal O_x \oplus
\mathcal O_x$ for some point $x \in Z'$.
\end{lemma}
\begin{proof}
Our lemma follows immediately from the following claim.

\medskip\noindent
{\bf Claim.}
{\it Let $Q$ be a length-$2$ torsion sheaf supported on
at most two points of $X$. Then, either $Q \cong
\mathcal O_{\xi}$ for some $\xi\in X^{[2]}$, or $Q \cong
\mathcal O_x \oplus \mathcal O_x$ for some $x \in X$.}

\medskip\noindent
{\it Proof.}
If $\Supp(Q) = \{x_1, x_2\}$ with $x_1 \ne x_2$, then $Q \cong
\mathcal O_{x_1} \oplus \mathcal O_{x_2}$ and we are done.
In the following, we assume that $\Supp(Q) = \{x\}$
for some $x \in X$. We further assume that
$X=\text{Spec}(A)$ is affine and $Q$ is
an $A$-module by abuse of notation.

Take a nonzero $v \in Q$ and define an $A$-module homomorphism
$\varphi\colon A \to Q$ by sending $1$ to $v$. For the ideal
$J=\ker\varphi$, the induced homomorphism $\overline\varphi\colon
A/J\to Q$ is injective. If $\overline\varphi$ is also surjective,
then $Q \cong \mathcal O_{\xi}$ for some $\xi\in X^{[2]}$. If
$\overline\varphi$ is not surjective, then we have an exact
sequence
\begin{eqnarray*}
0\rightarrow A/J \rightarrow Q\rightarrow Q^{\prime} \rightarrow 0
\end{eqnarray*}
where $Q^{\prime}$ and $A/J$ must be of length one. Therefore
$A/J\cong \mathcal O_x$ and $Q^{\prime}\cong \mathcal O_x$.

It follows that the minimal number of generators of $\mathcal Q$ is at
most two.

If $Q$ is generated by a single element, say $v_0$, then by
replacing the above $v \in Q$ by $v_0 \in Q$, we conclude that
$Q \cong \mathcal O_{\xi_0}$ for some $\xi_0 \in X^{[2]}$.

Now we are left with the case when the minimal number of
generators for $Q$ is two. Assume that $v_1$ and $v_2$ generate $Q$.
We define two homomorphisms:
\begin{eqnarray*}
&&\varphi_1\colon A\rightarrow Q,\quad \varphi_1(1)=v_1;\\
&&\varphi_2\colon A\rightarrow Q,\quad \varphi_2(1)=v_2.
\end{eqnarray*}
Then $\varphi_1$ and $\varphi_2$ induce injective homomorphisms:
\begin{eqnarray*}
\overline\varphi_1\colon A/J_1\rightarrow Q
\qquad \hbox{and }\qquad
\overline\varphi_2\colon A/J_2\rightarrow Q
\end{eqnarray*}
respectively. Note that both $A/J_1$ and $A/J_2$ must be of length-$1$.
So $A/J_1 \cong \mathcal O_x$ and $A/J_2 \cong \mathcal O_x$.
Define $\psi \colon A/J_1\oplus A/J_2\rightarrow Q$ by putting
\begin{eqnarray*}
\psi(a, b)=av_1+bv_2.
\end{eqnarray*}
Then $\psi$ is surjective. Since both $Q$ and $A/J_1\oplus A/J_2$
have length two, $\psi$ must be an isomorphism.
Therefore, we see that $Q \cong \mathcal O_x \oplus \mathcal O_x$.
\qed
\end{proof}

Conversely, we can show that if $[Z'] \in \mathfrak I_{1-g}(X, \beta)$
and $Z'$ is smooth at a point $x \in X$, then both types of $Q$
in Lemma~\ref{Q_str} can occur in (\ref{exact2}).

\begin{proposition} \label{xx}
Let $[Z'] \in \mathfrak I_{1-g}(X, \beta)$ and $Z'$ be smooth at
a point $x \in X$. Then,

{\rm (i)} there exists a unique $[Z] \in
\mathfrak I_{3-g}(X, \beta)$ sitting in the exact sequence:
\begin{eqnarray}   \label{xx.1}
0 \to I_Z \to I_{Z'} \to \mathcal O_x \oplus \mathcal O_x \to 0;
\end{eqnarray}

{\rm (ii)} $\dim \Hom(I_Z, \mathcal O_Z)
= 10 + \dim \Hom(I_{Z'}, \mathcal O_{Z'})$.

{\rm (iii)} the moduli space $\mathfrak I_{3-g}(X, \beta)$ is
not smooth at $[Z]$.
\end{proposition}
\begin{proof}
(i) It suffices to show the existence and uniqueness of $Z$ in
an analytic neighborhood $U_x$ of $x$. Let $w_1, w_2, w_3$ be
the coordinates of $U_x$ centered at $x$ such that $Z'$ is
given by $w_1 = w_2 = 0$. Then $I_{Z'} = (w_1, w_2) \subset
\C [w_1, w_2, w_3]$. Note that
\begin{eqnarray}   \label{xx.2.1}
(w_1, w_2) \otimes {\C [w_1, w_2, w_3] \over (w_1, w_2, w_3)}
\cong {(w_1, w_2) \over (w_1, w_2) \cdot (w_1, w_2, w_3)}
\cong \C \oplus \C.
\end{eqnarray}
Hence tensoring $\C [w_1, w_2, w_3] \to
{\C [w_1, w_2, w_3]/(w_1, w_2, w_3)} \to 0$ by $I_{Z'}$ yields
\begin{eqnarray}   \label{xx.3}
I_{Z'} \to \mathcal O_x \oplus \mathcal O_x \to 0.
\end{eqnarray}
The kernel of the surjection $I_{Z'} \to \mathcal O_x \oplus
\mathcal O_x$ defines an element $[Z] \in \mathfrak I_{3-g}(X, \beta)$
satisfying the exact sequence (\ref{xx.1}).
Note that in $U_x$, we have
\begin{eqnarray}  \label{xx.4}
I_Z = (w_1, w_2, w_3) \cdot (w_1, w_2)
    = (w_1^2, w_1w_2, w_2^2, w_1w_3, w_2w_3).
\end{eqnarray}

For the uniqueness of $Z$, note that specifying a surjection
$$I_{Z'} \to \mathcal O_x \oplus \mathcal O_x \to 0$$
is equivalent to specifying a surjection $I_{Z'} \otimes \mathcal O_x
\to \mathcal O_x \oplus \mathcal O_x \to 0$, i.e.,
$$\mathcal O_x \oplus \mathcal O_x \to
\mathcal O_x \oplus \mathcal O_x \to 0$$
in view of (\ref{xx.2.1}). Now the surjection
$\mathcal O_x \oplus \mathcal O_x \to
\mathcal O_x \oplus \mathcal O_x$ must be an isomorphism.
Therefore, there is only one quotient class
$[I_{Z'} \to O_x \oplus \mathcal O_x]$ up to isomorphisms of
$\mathcal O_x \oplus \mathcal O_x$.
This proves the uniqueness of $[Z] \in
\mathfrak I_{3-g}(X, \beta)$ satisfying (\ref{xx.1}).

(ii) We cover $X$ by the (analytic) open subsets $U_x$ and
$X - \{x\}$. Regard $\Hom(I_{Z'}, \mathcal O_{Z'})$ as obtained
from $\Hom(I_{Z'}|_{U_x}, \mathcal O_{Z'}|_{U_x})$ and
\begin{eqnarray}  \label{xx.5}
\Hom \big (I_{Z'}|_{X - \{x\}}, \mathcal O_{Z'}|_{X - \{x\}}
 \big )
\end{eqnarray}
by gluing along $U_x \cap (X - \{x\})$. Since $I_{Z'}|_{U_x}
= (w_1, w_2)$ and
$$\mathcal O_{Z'}|_{U_x} = {\C [w_1, w_2, w_3] \over (w_1, w_2)}
\cong \C [w_3],$$
we see that the homomorphisms $f \in \Hom(I_{Z'}|_{U_x},
\mathcal O_{Z'}|_{U_x})$ are of the form
\begin{eqnarray}  \label{xx.6}
f(w_1) = \sum_{i=0}^{+\infty} a_{1i} w_3^i, \qquad
f(w_2) = \sum_{i=0}^{+\infty} a_{2i} w_3^i
\end{eqnarray}
where $a_{1i},\ a_{2i}$ are independent complex parameters.

Similarly, $\Hom(I_{Z}, \mathcal O_{Z})$ is obtained
from $\Hom(I_{Z}|_{U_x}, \mathcal O_{Z}|_{U_x})$ and
\begin{eqnarray}  \label{xx.7}
\Hom \big (I_{Z}|_{X - \{x\}}, \mathcal O_{Z}|_{X - \{x\}}
 \big )
\end{eqnarray}
by gluing along $U_x \cap (X - \{x\})$. Note that (\ref{xx.5}) and
(\ref{xx.7}) are identical.

Next, we compare
$\Hom(I_{Z}|_{U_x}, \mathcal O_{Z}|_{U_x})$ and
$\Hom(I_{Z'}|_{U_x}, \mathcal O_{Z'}|_{U_x})$.
By (\ref{xx.4}), we check that
the homomorphisms $h \in \Hom(I_{Z}|_{U_x}, \mathcal O_{Z}|_{U_x})$
are of the form
\begin{eqnarray}
h(w_1^2)  &=& b_{11} w_1 + b_{12} w_2, \nonumber \\
h(w_1w_2) &=& b_{21} w_1 + b_{22} w_2, \nonumber \\
h(w_2^2)  &=& b_{31} w_1 + b_{32} w_2, \nonumber \\
h(w_1w_3) &=& b_{41} w_1 + b_{42} w_2 +
              \sum_{i=1}^{+\infty} a_{1i}' w_3^i, \label{xx.8.1} \\
h(w_2w_3) &=& b_{51} w_1 + b_{52} w_2 +
              \sum_{i=1}^{+\infty} a_{2i}' w_3^i    \label{xx.8.2}
\end{eqnarray}
where $b_{ij},\ a_{ij}'$ are independent
complex parameters. By (\ref{xx.4}) again,
\begin{eqnarray}  \label{xx.9}
I_{Z}|_{U_x \cap (X - \{x\})} = (w_1, w_2).
\end{eqnarray}
So $w_1 = w_2 = 0$ in $\mathcal O_{Z}|_{U_x \cap (X - \{x\})}$,
and the restrictions of the homomorphisms $h \in \Hom(I_{Z}|_{U_x},
\mathcal O_{Z}|_{U_x})$ to $U_x \cap (X - \{x\})$ are of the form
\begin{eqnarray*}
h(w_1w_3) = \sum_{i=1}^{+\infty} a_{1i}' w_3^i,  \qquad
h(w_2w_3) = \sum_{i=1}^{+\infty} a_{2i}' w_3^i.
\end{eqnarray*}
Combining this with (\ref{xx.9}), we see that
the restrictions of the homomorphisms $h \in \Hom(I_{Z}|_{U_x},
\mathcal O_{Z}|_{U_x})$ to $U_x \cap (X - \{x\})$ are of the form
\begin{eqnarray}  \label{xx.10}
h(w_1) = \sum_{i=0}^{+\infty} a_{1i}' w_3^i,  \qquad
h(w_2) = \sum_{i=0}^{+\infty} a_{2i}' w_3^i
\end{eqnarray}
which is precisely of the form (\ref{xx.6}). Therefore, we conclude that
\begin{eqnarray*}
   \dim \Hom(I_Z, \mathcal O_Z)
&=&\#\{b_{ij}|1 \le i \le 5, 1 \le j \le 2\}
     + \dim \Hom(I_{Z'}, \mathcal O_{Z'})          \\
&=&10 + \dim \Hom(I_{Z'}, \mathcal O_{Z'}).
\end{eqnarray*}

(iii) Let $\mathfrak I_{1-g}$ be an irreducible component of
the moduli space $\mathfrak I_{1-g}(X, \beta)$ containing the point $[Z']$.
Let $\mathfrak U$ be the open subset of the moduli space
$\mathfrak I_{3-g}(X, \beta)$ consisting of all the closed subschemes of
$X$ of the form:
$$C \cup \{x_1, x_2\}$$
where $[C] \in \mathfrak I_{1-g}$ and $x_1, x_2$ are distinct points
not contained in the curve $C$. Then $[Z] \in \overline{\mathfrak U}$
since $Z$ is the flat limit as $t$ approaches 0 of the subschemes
\begin{eqnarray}   \label{flat_lim}
Z' \cup \{x_1(t), \,\, x_2(t)\}
\end{eqnarray}
where $x_1(t) = (t,0,0) \in U_x \subset X$ and
$x_2(t) = (0,t,0) \in U_x \subset X$.
Hence $\overline{\mathfrak U}$ is an irreducible component
of the moduli space $\mathfrak I_{3-g}(X, \beta)$
containing the point $[Z]$. By (\ref{Zar_tang}),
the Zariski tangent space to $\mathfrak I_{3-g}(X, \beta)$ at $[Z]$
is $\text{Hom}(I_{Z}, \mathcal O_{Z})$,
and the Zariski tangent space to $\mathfrak I_{1-g}(X, \beta)$
at $[Z']$ is $\text{Hom}(I_{Z'}, \mathcal O_{Z'})$. Since
\begin{eqnarray*}
   \dim \Hom(I_Z, \mathcal O_Z)
&=  &10 + \dim \Hom(I_{Z'}, \mathcal O_{Z'})    \\
&\ge&10 + \dim \,\, \mathfrak I_{1-g}           \\
&=  &4 + \dim \,\, \overline{\mathfrak U},
\end{eqnarray*}
we conclude that the moduli space $\mathfrak I_{3-g}(X, \beta)$
is not smooth at $[Z]$.
\end{proof}

Similarly, assuming the same conditions as in Proposition~\ref{xx}
and assuming $\xi \in X^{[2]}$ with $\Supp(\xi) = \{x\}$,
we can prove that if $\xi$ points to the tangent direction
of $Z'$ at $x$, then the set of surjections $I_{Z'} \to
\mathcal O_\xi \to 0$ up to isomorphisms has dimension $2$.
If $\xi$ is transverse to $Z'$ at $x$, then the set of surjections
$I_{Z'} \to \mathcal O_\xi \to 0$ up to isomorphisms has dimension $1$.
However, it is not clear how to globalize these local data into
a global description of $\mathfrak I_{3-g}(X, \beta)$
in terms of $\mathfrak I_{1-g}(X, \beta)$.

\begin{remark} \label{dt_euler}

\

\smallskip\noindent
(i) By Proposition~\ref{xx}~(iii), the moduli space $\mathfrak
I_{3-g}(X, \beta)$ is not smooth in general. Hence there is no
guarantee that the corresponding Donaldson-Thomas invariant (when
$K_X = 0$) is equal to the topological Euler number of $\mathfrak
I_{3-g}(X, \beta)$ up to sign.  Note however that the equality
does occur in some important cases with a singular moduli space.
This is most notably the case for the degree-$0$ Donaldson-Thomas
invariants. These degree-$0$ invariants have been conjectured in
\cite{MNOP1, MNOP2}, and computed by Jun Li \cite{Li}.

\smallskip\noindent
(ii) Under Assumption~\ref{assumption}, the topological Euler
number of the moduli space
$$\mathfrak I_{1, 2} = \mathfrak I_{3-g}(X, \beta)$$
(see (\ref{Imn})) has been computed in \cite{LQ2} by using virtual
Hodge polynomials.

\end{remark}

\end{document}